\documentclass [11pt]{article}
\usepackage{latexsym}

\newcommand{\qed}{\mbox{$\Diamond$}\vspace{\baselineskip}}

\newtheorem{theorem}{Theorem}

\newtheorem{lemma}{Lemma}

\newenvironment{proof}{\noindent {\bf Proof:}}{{\qed}}

\newcommand{\vanish}[1]{}
 
\begin{document}
\title{Packing Ferrers Shapes}

\author{Noga Alon
\thanks{School of Mathematics,
Institute for Advanced
Study, Princeton, NJ 08540, and
Raymond and Beverly Sackler Faculty of Exact Sciences,
Tel Aviv University, Tel Aviv, Israel. Email: noga@math.tau.ac.il.
Research supported in part by a Sloan Foundation grant 96-6-2, by a
State of New Jersey grant and by a USA-Israel BSF grant.}
\and Mikl\'os B\'ona \thanks{
School of Mathematics, Institute for Advanced Study, Princeton, NJ 08540. 
Supported by Trustee Ladislaus von Hoffmann, the Arcana Foundation.
Email: bona@math.ias.edu.}
\and Joel Spencer \thanks{
School of Mathematics, Institute for Advanced Study, Princeton, NJ 08540, and 
Department of Mathematics and Computer Science, Courant
Institute, NYU, NY. 
Email:spencer@cs.nyu.edu.
Research supported in part by a USA-Israel BSF Grant. }}

\maketitle

\begin{abstract}
Answering a question of Wilf, we show that if $n$ is sufficiently large,
then one cannot cover an $n \times  p(n)$ rectangle using 
each of the $p(n)$ distinct Ferrers
shapes of size $n$ exactly once. Moreover, the maximum number of pairwise
distinct, non-overlapping Ferrers shapes that can be packed in such
a rectangle is only $\Theta( p(n)/ \log n).$
\end{abstract}

\section{Introduction}
A {\em partition} $p$ of a positive integer $n$ is an  array $p=(x_1,x_2,\cdots
,x_k)$ of positive integers so that 
$x_1\geq x_2 \geq \cdots \geq x_k$ and
 $n=\sum_{i=1}^k x_i$. The $x_i$
are called the {\em parts} of $p$. The total number of distinct partitions
of $n$ is denoted by $p(n)$. A {\em Ferrers shape} of a partition 
$p=(x_1,x_2,\cdots,x_k)$ is a set of $n$ square boxes with sides parallel 
to the
coordinate axes so that in the $i$th row we have $x_i$ boxes and all rows 
start at the same vertical line. The Ferrers
shape of the partition $p=(4,2,1)$ is shown in Figure 1. Clearly, there is
an obvious bijection between partitions of $n$ and Ferrers shapes of size $n$.

If we reflect a Ferrers shape of a partition $p$ with respect 
to its main diagonal,
we get another shape, representing the {\em conjugate} partition of $p$.
Thus, in our example, the conjugate of (4,2,1) is (3,2,1,1).

\begin{center}

 \setlength{\unitlength}{0.00030000in}%
\begingroup\makeatletter\ifx\SetFigFont\undefined
\def\x#1#2#3#4#5#6#7\relax{\def\x{#1#2#3#4#5#6}}%
\expandafter\x\fmtname xxxxxx\relax \def\y{splain}%
\ifx\x\y   
\gdef\SetFigFont#1#2#3{%
  \ifnum #1<17\tiny\else \ifnum #1<20\small\else
  \ifnum #1<24\normalsize\else \ifnum #1<29\large\else
  \ifnum #1<34\Large\else \ifnum #1<41\LARGE\else
     \huge\fi\fi\fi\fi\fi\fi
  \csname #3\endcsname}%
\else
\gdef\SetFigFont#1#2#3{\begingroup
  \count@#1\relax \ifnum 25<\count@\count@25\fi
  \def\x{\endgroup\@setsize\SetFigFont{#2pt}}%
  \expandafter\x
    \csname \romannumeral\the\count@ pt\expandafter\endcsname
    \csname @\romannumeral\the\count@ pt\endcsname
  \csname #3\endcsname}%
\fi
\fi\endgroup
\begin{picture}(5124,4074)(2389,-5323)
\thicklines
\put(2401,-2611){\framebox(1275,1350){}}
\put(3676,-2611){\framebox(1275,1350){}}
\put(4951,-2611){\framebox(1275,1350){}}
\put(6226,-2611){\framebox(1275,1350){}}
\put(2401,-3961){\framebox(1275,1350){}}
\put(3676,-3961){\framebox(1275,1350){}}
\put(2401,-5311){\framebox(1275,1350){}}
\end{picture}
\end{center}
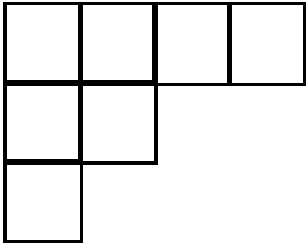
\centerline{Figure 1: The Ferrers shape of (4,2,1)}

Recently, Herb Wilf \cite{wilf} has asked the following intriguing question:
consider all distinct Ferrers shapes consisting of $n$ boxes.
 Is it true that  for $n$ sufficiently large, one can always
 tile a rectangle of side lengths $p(n)$ and $n$ using each of these shapes
 exactly once ? Obviously, in such a tiling, if one exists, the shapes cannot
overlap each other. For small values of $n$, one gets mixed answers:
for $n=1,2,4$ such a tiling exists, however, for $n=3$ there is no such 
tiling.

In this short paper we  answer Wilf's question in the negative by showing 
that for $n$ sufficiently large no such tiling exists. In fact, we prove
the following stronger statement.
\begin{theorem} 
\label{t1}
If $n$ is sufficiently large, then one cannot cover an 
$n \times  p(n)$ rectangle by using each of the $p(n)$ 
distinct Ferrers
shapes of size $n$ exactly once. Moreover, the maximum fraction
of the area of this rectangle that can be covered by non-overlapping
distinct Ferrers shapes of size $n$ is   at most
$\frac{c}{\log n}$, for some absolute constant $c$.
\end{theorem}

The $c/\log n$ upper bound 
is tight, up to the constant $c$, and shows that
as $n$ grows we cannot even cover
a fixed fraction of the area by non-overlapping distinct shapes.

To prove the result, we use some geometric properties that are shared
by the vast majority of the Ferrers shape of size $n$ and imply that
these shapes cannot be packed in an efficient way.
The geometric properties we need can probably be derived from 
the extensive available information
on the typical form of a Ferrers shape, given, for example, in
 \cite{pittel}, and partly also in several earlier papers including
 \cite{szt1},  \cite{szt2},
 \cite{szt3}.
However, 
since we need an explicit estimate of the number of shapes that
do not satisfy the required properties, and in order to make the paper
self-contained, we prefer to derive all of them directly from the
Hardy-Ramanujan asymptotic formula for $p(n)$.
This is done in the next section. In Section 3, we apply the geometric
properties to prove our main result. Throughout the paper we assume,
whenever this is needed, that the size $n$ of the Ferrers shapes
considered is sufficiently large.

\section{Some geometric properties of typical Ferrers shapes }

The corner of the 
first row and first column of a Ferrers shape will be called
the {\em apex} of that shape.
In this section we prove
some asymptotic geometric properties of Ferrers shapes of size $n$. 
Our basic tool is the well-know 
 Hardy-Ramanujan asymptotic formula for the number of shapes of size
$n$, see, e.g., \cite{andrews}.
It asserts that
\begin{equation} \label{rama}
p(n) =(1+o(1)) \frac{e^{C\sqrt{n}}}{4n\sqrt{3}}, \end{equation}
where $C=\pi \sqrt{\frac{2}{3}}$, and the $o(1)$-term tends to
$0$ as $n$ tends to infinity.

\begin{lemma} 
\label{l1} 
Let $x_1 \geq x_2 \geq \ldots \geq x_k$ denote the parts of
the partition corresponding to the Ferrers shape $F$, and let
$y_1 \geq y_2 \ldots \geq y_s$ denote  the parts of the 
conjugate partition.
The following hold for  all but $p(n)/n^{0.48}$ 
Ferrers shapes $F$ of  size $n$.
\begin{itemize}
\item{I.} There exist absolute constants $c_1$ and $c_2$ so that for $n$ 
sufficiently large, we
 have $c_1 \sqrt n \log n < x_1 < c_2 \sqrt n \log n,$ and also,
$c_1 \sqrt n \log n < y_1 < c_2 \sqrt n \log n $.
\item{II.} 
For every $\epsilon >0$ there exists a $\delta>0$ such that the
following holds.
Let $Q$ be a square of side length $\epsilon \sqrt n \log n$ whose
 apex is that of $F$ and whose sides are adjusted to the border of $F$,
(see Figure 2). Then 
$F - Q$ contains at most $n^{1-\delta}$ boxes.
\item{III.} There exists $c_3>0$ so that
 $F$ has at least $c_3\sqrt n$ parts of size at least $c_3\sqrt n$ each.
\end{itemize}
\end{lemma}
\begin{proof}
\begin{itemize}
\item{I.} 
The inequalities for $x_1$ and $y_1$ are clearly equivalent by taking
conjugates. Thus it suffices to prove the statements for $x_1$.
\begin{itemize}
\item We first bound the number of partitions with a large value of
$x_1$. To do so,
consider a partition with 
$x_1 > c_2 \sqrt n \log n$ and remove its first part. This way
we get a partition of an integer smaller than $n-c_2 \sqrt n \log n$.
There are (less than) $n$ ways to choose the first part, and by the
obvious monotonicity of $p(n)$ there are at most 
$p(n-c_2 \sqrt{n} \log n)$ ways to choose the remaining partition.
Therefore, the number of partitions of $n$ with $x_1 > c_2 \sqrt n \log n$
does not exceed
$n\cdot p(n-c_2 \sqrt{n} \log n)$.
By (\ref{rama}) this yields
\[n\cdot p(n-c_2 \sqrt n \log n)\leq n 
\frac{e^{C\sqrt{n-c_2 \sqrt n \log n}}}{4n\sqrt 3}<n \frac{p(n)}{n^{0.5Cc_2}},
 \]
as $\sqrt{n-c_2 \sqrt n \log n}\leq \sqrt{n}-\frac{1}{2}c_2 \log n$. 
Since $C=\pi \sqrt {2/3}$ this implies that by choosing,
for example, 
$c_2$ to be $4$, we conclude that the number of partitions of $n$
that do not satisfy
$x_1 \leq c_2 \sqrt n \log n$ is smaller than $p(n)/n$.

\item Next we prove that for almost all partitions we also have  
$x_1 \leq  c_1 \sqrt n \log n$, for some positive constant $c_1$.
 Let $S$ be the set of partitions of $n$ violating this
constraint, and attach two additional 
parts $x_0$ and $x_{-1}$ in
all possible ways to all partitions in $S$ so that the following hold:
\[x_0+x_{-1}=3\cdot[ c_1 \sqrt n \log n],\]
and 
\[x_{-1}\geq x_{0}\geq c_1 \sqrt n \log n.\]
Let $S'$ be the set of partitions obtained this way. It then follows that
$x_{-1}$ and $x_{0}$ are the two largest parts in all partitions
in $S'$, and that $S'$ contains partitions of the integer 
$n+3\cdot[ c_1 \sqrt n \log n]$.

As $x_{-1}\geq x_0$, we must have $1.5\cdot c_1 \sqrt n \log n
\leq x_{-1} \leq 2\cdot c_1 \sqrt n \log n$, so we have 
 $0.5\cdot c_1 \sqrt n \log n$ choices for $x_{-1}$.
This implies 
\[|S'|= |S|\cdot 0.5\cdot c_1 \sqrt n \log n
 \leq p(n+3\cdot[ c_1 \sqrt n \log n]),\]
which yields
\[|S|=\frac{|S'|}{0.5\cdot c_1 \sqrt n \log n}
 \leq 
\frac{p(n+3\cdot[ c_1 \sqrt n \log n])}{0.5\cdot c_1 \sqrt n \log n}\leq
\frac{p(n)\cdot e^{1.5Cc_1 \log n}}{0.5\cdot c_1 \sqrt n \log n} \leq
\frac{p(n)\cdot n^{1.5Cc_1}}{0.5\cdot c_1 \sqrt n \log n}, \]
as $\sqrt{n+3c_1 \sqrt n \log n}< \sqrt n +1.5c_1\log n$. By 
choosing a sufficiently small $c_1$  
(e.g., $c_1=1/1000$), we see that  all but $p(n)/n^{0.49}$
 partitions satisfy $x_1 \geq c_1 \sqrt n \log n$.
(Note that it is not difficult to show that in fact for every fixed $r$
there is some $c_1>0$ so that for all but at most $p(n)/n^r$ partitions
of $n$, $x_1 \geq c_1 \sqrt n \log n$. This can be done by adding to
each partition that does not satisfy the above more than $2$ parts with
a prescribed sum, and by repeating the above argument. For our
purpose here, however, the above estimate suffices.)
\end{itemize}
\item{II.}  It suffices to prove that most partitions do not have
too many parts which are larger than $\epsilon \sqrt n \log n$,
as this would imply, by taking the conjugate, 
that most partitions also do not have
too many columns larger than $\epsilon \sqrt n \log n$, either, 
yielding that the total area outside the square $Q$ is small. 
More precisely, 
 we prove that for every $\epsilon>0$ there exists a $\delta>0$ so that
for all but at most $p(n)/n$ of the partitions of $n$,
$|\{i |x_i >\epsilon \sqrt n \log n\}|<n^{\frac{1}{2}-\delta}$.
\begin{center}
\setlength{\unitlength}{0.00027500in}%
\begingroup\makeatletter\ifx\SetFigFont\undefined
\def\x#1#2#3#4#5#6#7\relax{\def\x{#1#2#3#4#5#6}}%
\expandafter\x\fmtname xxxxxx\relax \def\y{splain}%
\ifx\x\y   
\gdef\SetFigFont#1#2#3{%
  \ifnum #1<17\tiny\else \ifnum #1<20\small\else
  \ifnum #1<24\normalsize\else \ifnum #1<29\large\else
  \ifnum #1<34\Large\else \ifnum #1<41\LARGE\else
     \huge\fi\fi\fi\fi\fi\fi
  \csname #3\endcsname}%
\else
\gdef\SetFigFont#1#2#3{\begingroup
  \count@#1\relax \ifnum 25<\count@\count@25\fi
  \def\x{\endgroup\@setsize\SetFigFont{#2pt}}%
  \expandafter\x
    \csname \romannumeral\the\count@ pt\expandafter\endcsname
    \csname @\romannumeral\the\count@ pt\endcsname
  \csname #3\endcsname}%
\fi
\fi\endgroup
\begin{picture}(4524,4824)(3664,-6073)
\thicklines
\put(3676,-5236){\framebox(3900,3975){}}
\put(3676,-5236){\line( 0,-1){825}}
\put(3676,-6061){\line( 1, 0){750}}
\put(4426,-6061){\line( 0, 1){1425}}
\put(4426,-4636){\line( 1, 0){900}}
\put(5326,-4636){\line( 0, 1){975}}
\put(5326,-3661){\line( 1, 0){1350}}
\put(6676,-3661){\line( 0, 1){1575}}
\put(6676,-2086){\line( 1, 0){900}}
\put(7576,-1261){\line( 1, 0){600}}
\put(8176,-1261){\line( 0,-1){825}}
\put(7576,-2086){\line( 1, 0){600}}
\put(7276,-4936){\makebox(0,0)[lb]{\smash{\SetFigFont{7}{8.4}{rm}Q}}}
\put(4951,-4411){\makebox(0,0)[lb]{\smash{\SetFigFont{7}{8.4}{rm}F}}}
\end{picture}

\end{center}
\centerline{Figure 2: $Q$ contains most of $F$}

Let $S=\{p~:~\mbox{ $p=(x_1,\ldots ,x_k)$ is a partition of $n$ so 
that } |\{i |x_i >\epsilon \sqrt n \log n
\}| \geq n^{1/2-\delta} \}$. As we have already proved statement I of the 
lemma, we may and will assume that all parts of the members of $S$ 
are smaller than $c_2 \sqrt n \log n$, as this holds for all of them but 
at most $p(n)/n$. Note, also, that if the maximum row and column length
of a partition is at most $c_2 \sqrt n \log n$, and it has at most
$n^{1/2-\delta}$ rows and at most that many columns of length
exceeding $\epsilon \sqrt n \log n$,  then the total area of its Ferrers
shape that is not covered by the square $Q$ is at most
$2c_2n^{1-\delta} \log n<n^{1-\delta'}$ for, say, $\delta'=\delta/2$,
as needed.

To obtain a partition in $S$, 
we first need to choose $n^{\frac{1}{2}-\delta}$ parts,
each larger than $\epsilon \sqrt n \log n$ and smaller than 
$c_2 \sqrt n \log n$. This can be done in at most
 ${{ c_2 \sqrt n \log n
+ n^{\frac{1}{2}-\delta}\choose n^{\frac{1}{2}-\delta} }}$ different ways. 
Next, we choose
a partition of the remaining integer. This can be done  in at most 
$p(n-\epsilon n^{1-\delta}\log n)$ ways. Therefore,
\begin{equation} \label{square} |S|\leq {{ c_2 \sqrt n \log n
+ n^{\frac{1}{2}-\delta}\choose n^{\frac{1}{2}-\delta} }} \cdot 
p(n-\epsilon n^{1-\delta}\log n).\end{equation}
By Stirling's formula 
\[{{a \choose b}} \leq (\frac{ea}{b})^b.\]
Thus inequality (\ref{square}) implies
\[|S|<(2ec_2n^{\delta}\log n )^{n^{0.5-\delta}} \cdot p(n)e^{-0.5C\epsilon 
n^{0.5-\delta}\log n}= p(n)(\frac{2ec_2n^{\delta}\log n }{n^{0.5C\epsilon}})^
{n^{0.5-\delta}}.\]
 As $C>1$, choosing, e.g., $\delta= \epsilon /2$ gives the desired estimate.
\item{III.} Let $S=\{~p:~ p=(x_1,\ldots, x_k) \mbox { is a 
partition of $n$ so that } |\{i|x_i\geq c_3\sqrt{n}\} |
<c_3 \sqrt n \}$. 
Let $\cal F$ be a family of at least, say, $2^{0.1c_3 \sqrt n}$
subsets of a set of
cardinality $10c_3 \sqrt n$ in which the Hamming distance between any
two subsets is larger than $2c_3 \sqrt n$. (The existence of such  a family
is easy and follows from the Gilbert-Varshamov bound, see, e.g., 
\cite{MS}.)
Let the
underlying set of $\cal F$ be the set $\{c_3 \sqrt n+1, 
c_3 \sqrt n +2, \cdots ,  11c_3\sqrt n\}$.

Define $S'=\{P\cup F | P\in S, 
F\in \cal F\}$. It is not too difficult to check that all
$|S|\cdot |\cal F|$ partitions in $S'$ are pairwise distinct;
indeed, if two such unions have the same $P$ or the same $F$ then they
clearly differ. On the other hand, for distinct $P,P'$ in $S$ and distinct
$F,F'$ in ${\cal F}$, $P \cup F$ and $P' \cup F'$ do not have the same sets of
parts of size bigger than $c_3 \sqrt n$, by the definition of $S$ and the
choice of ${\cal F}$.
Therefore, we have
\[|S'|=|S|\cdot |\cal F| \]
which yields
\begin{equation}\label{third} 
|S|\leq \frac{|S'|}{2^{0.1c_3 \sqrt n}}. \end{equation}
All the elements of $S'$
are partitions of integers not larger than $n+110c_3^2n$.
As these are all distinct it follows that
\[|S'|\leq \sum_{k\leq n+110c_3^2n}p(k)\leq 
e^{C\sqrt{n(1+110c_3^2})}\leq e^{C\sqrt{n}+55C\cdot c_3^2\sqrt{n}},\]
and therefore, by inequality (\ref{third}),
\[|S|\leq \frac{e^{C\sqrt{n}}e^{55C\cdot c_3^2\sqrt{n}}}{2^{0.1c_3 \sqrt n}}\]
which completes the proof, as by choosing, say, $c_3=0.001$, 
the second term of the
numerator becomes much smaller than the denominator.
\end{itemize}
This completes the proof of Lemma \ref{l1}.
\end{proof}

\section{The proof of the main result}

In this section we prove Theorem \ref{t1}. Consider a collection
of non-overlapping pairwise distinct Ferrers shapes of size $n$
contained in a given $n$ by $p(n)$ rectangle $R$. First,
omit  all shapes that do not satisfy the three conditions in
Lemma \ref{l1}, and call the remaining collection of shapes the
{\em reduced collection}. By the assertion of the lemma at most $p(n)/n^{0.48}
<p(n)/\log n$ shapes are deleted in this step, and we are now left
only with shapes satisfying all three 
conditions of the lemma. 

Let $A\subset R$ be a square of side length $0.8c_1 \sqrt n \log n$ where
$c_1$ is a constant satisfying condition I of lemma \ref{l1}.

\begin{lemma}
\label{l2}
There exists an absolute constant $\gamma$ so that the total area of $A$ 
covered by our reduced collection is at most $\gamma n\log n$.
\end{lemma}
\begin{proof} First note that we can suppose all of our shapes have their
apexes in their lower left corner. Indeed, by symmetry, we can assume at least
1/4 of all shapes are like that, and losing a constant factor does not
change the validity of the lemma. 
Let us view $A$  as a center of a larger square $B$ of 
$c_1\sqrt n \log n$ by 
$ c_1\sqrt n \log n$, as shown in Figure 3. 

\begin{center}
\setlength{\unitlength}{0.00030000in}%
\begingroup\makeatletter\ifx\SetFigFont\undefined
\def\x#1#2#3#4#5#6#7\relax{\def\x{#1#2#3#4#5#6}}%
\expandafter\x\fmtname xxxxxx\relax \def\y{splain}%
\ifx\x\y   
\gdef\SetFigFont#1#2#3{%
  \ifnum #1<17\tiny\else \ifnum #1<20\small\else
  \ifnum #1<24\normalsize\else \ifnum #1<29\large\else
  \ifnum #1<34\Large\else \ifnum #1<41\LARGE\else
     \huge\fi\fi\fi\fi\fi\fi
  \csname #3\endcsname}%
\else
\gdef\SetFigFont#1#2#3{\begingroup
  \count@#1\relax \ifnum 25<\count@\count@25\fi
  \def\x{\endgroup\@setsize\SetFigFont{#2pt}}%
  \expandafter\x
    \csname \romannumeral\the\count@ pt\expandafter\endcsname
    \csname @\romannumeral\the\count@ pt\endcsname
  \csname #3\endcsname}%
\fi
\fi\endgroup
\begin{picture}(6099,6255)(1789,-6973)
\thicklines
\put(1801,-6961){\framebox(6075,6000){}}
\put(2476,-6361){\framebox(4800,4800){}}
\put(2176,-886){\makebox(0,0)[lb]{\smash{\SetFigFont{12}{14.4}{bf}B}}}
\put(2551,-1486){\makebox(0,0)[lb]{\smash{\SetFigFont{12}{14.4}{bf}A}}}
\end{picture}
\end{center}
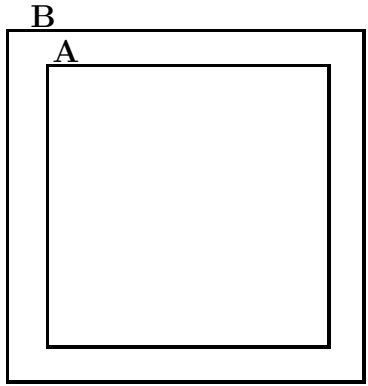
\centerline{Figure 3}

First, note that  shapes 
intersecting $A$ but having
their apexes outside $B$ cover very little of $A$. Indeed, by 
Lemma \ref{l1}, I all their
 apexes lie in a square of $(c_2 +c_1)\sqrt n \log n \times
 (c_2+c_1) \sqrt n \log n$, hence there are only $(2c_2+c_1)^2
\cdot \log^2 n$ of 
them, as they are
pairwise disjoint, the area of each of them is $n$,
and they all lie completely in a square of side length 
$(2c_2+c_1)\sqrt n \log n$. 
By Lemma \ref{l1}, II, each of them covers at most an area 
$n^{1-\delta}$ in $A$, and hence their total contribution is 
at most $O(n^{1-\delta} \log^2 n) < O(n)$.

Now consider the shapes whose apexes lie in $B$. Let  
$F_1,F_2$ be two such shapes, let $(x_1,y_1)$
be the apex of $F_1$, and 
 let $(x_2,y_2)$
be the apex of $F_2$,and suppose $x_1 \leq x_2$. Then by Lemma \ref{l1},I, 
$y_2$ cannot satisfy $y_2 \leq y_1$ as otherwise the shapes $F_1$ and 
$F_2$ would intersect. 
Therefore, $y_2\geq y_1$ must hold. 
 It follows that if we sort these shapes according to their increasing 
$x$ coordinates, the order we
obtain is also an increasing order of their $y$ coordinates. Moreover, by
Lemma \ref{l1}, III it follows
that $x_2+y_2 \geq x_1+y_1+c_3\sqrt n$, and a similar inequality holds
for each consecutive pair of 
shapes with apexes in $B$ in the above sorted order.
Hence, if $(x_i,y_i)$ is the apex of shape number $i$ in this
order,
$x_i+y_i\geq x_1+y_1+ (i-1) c_3\sqrt n$. Thus the total number of 
apexes in $B$ is at most $2 c_1\log n/c_3+1$, and therefore the 
total area of $A$ covered is not more than 
$\gamma n \log n$, for, say  $\gamma=3c_1/c_3$, completing the proof
of the lemma.
\end{proof} 

Theorem \ref{t1} is an immediate consequence of 
Lemma \ref{l2}, as one can essentially tile the $n$ by $p(n)$ rectangle
with squares of side length $0.8 c_1 \sqrt n \log n$. 
It is not difficult to  see that the assertion of the theorem is 
tight, up to the multiplicative constant $c$. Indeed, one can first omit
all shapes violating the condition in Lemma \ref{l1}, I, and then pack
the shapes along diagonals, where the apex of each shape touches
the furthest point on the main diagonal of the previous shape. 

\begin{center}
\setlength{\unitlength}{0.00033300in}%
\begingroup\makeatletter\ifx\SetFigFont\undefined
\def\x#1#2#3#4#5#6#7\relax{\def\x{#1#2#3#4#5#6}}%
\expandafter\x\fmtname xxxxxx\relax \def\y{splain}%
\ifx\x\y   
\gdef\SetFigFont#1#2#3{%
  \ifnum #1<17\tiny\else \ifnum #1<20\small\else
  \ifnum #1<24\normalsize\else \ifnum #1<29\large\else
  \ifnum #1<34\Large\else \ifnum #1<41\LARGE\else
     \huge\fi\fi\fi\fi\fi\fi
  \csname #3\endcsname}%
\else
\gdef\SetFigFont#1#2#3{\begingroup
  \count@#1\relax \ifnum 25<\count@\count@25\fi
  \def\x{\endgroup\@setsize\SetFigFont{#2pt}}%
  \expandafter\x
    \csname \romannumeral\the\count@ pt\expandafter\endcsname
    \csname @\romannumeral\the\count@ pt\endcsname
  \csname #3\endcsname}%
\fi
\fi\endgroup
\begin{picture}(7299,6549)(2464,-7273)
\thicklines
\put(2476,-736){\line( 0,-1){2400}}
\put(2476,-3136){\line( 1, 0){225}}
\put(2701,-3136){\line( 0, 1){1200}}
\put(2701,-1936){\line( 1, 0){225}}
\put(2926,-1936){\line( 0, 1){825}}
\put(2926,-1111){\line( 1, 0){675}}
\put(3601,-1111){\line( 0, 1){150}}
\put(3601,-961){\line( 1, 0){1200}}
\put(4801,-961){\line( 0, 1){225}}
\put(4801,-736){\line(-1, 0){2325}}
\put(2476,-736){\line( 0, 1){  0}}
\put(2926,-1111){\line( 0,-1){2400}}
\put(2926,-3511){\line( 1, 0){225}}
\put(3151,-3511){\line( 0, 1){1200}}
\put(3151,-2311){\line( 1, 0){225}}
\put(3376,-2311){\line( 0, 1){825}}
\put(3376,-1486){\line( 1, 0){675}}
\put(4051,-1486){\line( 0, 1){150}}
\put(4051,-1336){\line( 1, 0){1200}}
\put(5251,-1336){\line( 0, 1){225}}
\put(5251,-1111){\line(-1, 0){2325}}
\put(2926,-1111){\line( 0, 1){  0}}
\put(3376,-1486){\line( 0,-1){2400}}
\put(3376,-3886){\line( 1, 0){225}}
\put(3601,-3886){\line( 0, 1){1200}}
\put(3601,-2686){\line( 1, 0){225}}
\put(3826,-2686){\line( 0, 1){825}}
\put(3826,-1861){\line( 1, 0){675}}
\put(4501,-1861){\line( 0, 1){150}}
\put(4501,-1711){\line( 1, 0){1200}}
\put(5701,-1711){\line( 0, 1){225}}
\put(5701,-1486){\line(-1, 0){2325}}
\put(3376,-1486){\line( 0, 1){  0}}
\put(3826,-1861){\line( 0,-1){2400}}
\put(3826,-4261){\line( 1, 0){225}}
\put(4051,-4261){\line( 0, 1){1200}}
\put(4051,-3061){\line( 1, 0){225}}
\put(4276,-3061){\line( 0, 1){825}}
\put(4276,-2236){\line( 1, 0){675}}
\put(4951,-2236){\line( 0, 1){150}}
\put(4951,-2086){\line( 1, 0){1200}}
\put(6151,-2086){\line( 0, 1){225}}
\put(6151,-1861){\line(-1, 0){2325}}
\put(3826,-1861){\line( 0, 1){  0}}
\put(4276,-2236){\line( 0,-1){2400}}
\put(4276,-4636){\line( 1, 0){225}}
\put(4501,-4636){\line( 0, 1){1200}}
\put(4501,-3436){\line( 1, 0){225}}
\put(4726,-3436){\line( 0, 1){825}}
\put(4726,-2611){\line( 1, 0){675}}
\put(5401,-2611){\line( 0, 1){150}}
\put(5401,-2461){\line( 1, 0){1200}}
\put(6601,-2461){\line( 0, 1){225}}
\put(6601,-2236){\line(-1, 0){2325}}
\put(4276,-2236){\line( 0, 1){  0}}
\put(4726,-2611){\line( 0,-1){2400}}
\put(4726,-5011){\line( 1, 0){225}}
\put(4951,-5011){\line( 0, 1){1200}}
\put(4951,-3811){\line( 1, 0){225}}
\put(5176,-3811){\line( 0, 1){825}}
\put(5176,-2986){\line( 1, 0){675}}
\put(5851,-2986){\line( 0, 1){150}}
\put(5851,-2836){\line( 1, 0){1200}}
\put(7051,-2836){\line( 0, 1){225}}
\put(7051,-2611){\line(-1, 0){2325}}
\put(4726,-2611){\line( 0, 1){  0}}
\put(5176,-2986){\line( 0,-1){2400}}
\put(5176,-5386){\line( 1, 0){225}}
\put(5401,-5386){\line( 0, 1){1200}}
\put(5401,-4186){\line( 1, 0){225}}
\put(5626,-4186){\line( 0, 1){825}}
\put(5626,-3361){\line( 1, 0){675}}
\put(6301,-3361){\line( 0, 1){150}}
\put(6301,-3211){\line( 1, 0){1200}}
\put(7501,-3211){\line( 0, 1){225}}
\put(7501,-2986){\line(-1, 0){2325}}
\put(5176,-2986){\line( 0, 1){  0}}
\put(5626,-3361){\line( 0,-1){2400}}
\put(5626,-5761){\line( 1, 0){225}}
\put(5851,-5761){\line( 0, 1){1200}}
\put(5851,-4561){\line( 1, 0){225}}
\put(6076,-4561){\line( 0, 1){825}}
\put(6076,-3736){\line( 1, 0){675}}
\put(6751,-3736){\line( 0, 1){150}}
\put(6751,-3586){\line( 1, 0){1200}}
\put(7951,-3586){\line( 0, 1){225}}
\put(7951,-3361){\line(-1, 0){2325}}
\put(5626,-3361){\line( 0, 1){  0}}
\put(6076,-3736){\line( 0,-1){2400}}
\put(6076,-6136){\line( 1, 0){225}}
\put(6301,-6136){\line( 0, 1){1200}}
\put(6301,-4936){\line( 1, 0){225}}
\put(6526,-4936){\line( 0, 1){825}}
\put(6526,-4111){\line( 1, 0){675}}
\put(7201,-4111){\line( 0, 1){150}}
\put(7201,-3961){\line( 1, 0){1200}}
\put(8401,-3961){\line( 0, 1){225}}
\put(8401,-3736){\line(-1, 0){2325}}
\put(6076,-3736){\line( 0, 1){  0}}
\put(6526,-4111){\line( 0,-1){2400}}
\put(6526,-6511){\line( 1, 0){225}}
\put(6751,-6511){\line( 0, 1){1200}}
\put(6751,-5311){\line( 1, 0){225}}
\put(6976,-5311){\line( 0, 1){825}}
\put(6976,-4486){\line( 1, 0){675}}
\put(7651,-4486){\line( 0, 1){150}}
\put(7651,-4336){\line( 1, 0){1200}}
\put(8851,-4336){\line( 0, 1){225}}
\put(8851,-4111){\line(-1, 0){2325}}
\put(6526,-4111){\line( 0, 1){  0}}
\put(6976,-4486){\line( 0,-1){2400}}
\put(6976,-6886){\line( 1, 0){225}}
\put(7201,-6886){\line( 0, 1){1200}}
\put(7201,-5686){\line( 1, 0){225}}
\put(7426,-5686){\line( 0, 1){825}}
\put(7426,-4861){\line( 1, 0){675}}
\put(8101,-4861){\line( 0, 1){150}}
\put(8101,-4711){\line( 1, 0){1200}}
\put(9301,-4711){\line( 0, 1){225}}
\put(9301,-4486){\line(-1, 0){2325}}
\put(6976,-4486){\line( 0, 1){  0}}
\put(7426,-4861){\line( 0,-1){2400}}
\put(7426,-7261){\line( 1, 0){225}}
\put(7651,-7261){\line( 0, 1){1200}}
\put(7651,-6061){\line( 1, 0){225}}
\put(7876,-6061){\line( 0, 1){825}}
\put(7876,-5236){\line( 1, 0){675}}
\put(8551,-5236){\line( 0, 1){150}}
\put(8551,-5086){\line( 1, 0){1200}}
\put(9751,-5086){\line( 0, 1){225}}
\put(9751,-4861){\line(-1, 0){2325}}
\put(7426,-4861){\line( 0, 1){  0}}
\put(2476,-3136){\line( 6,-5){4950}}
\put(4801,-736){\line( 6,-5){4950}}
\put(2476,-7261){\framebox(7275,6525){}}
\put(2476,-5686){\line( 6,-5){1836.885}}
\put(7801,-736){\line( 5,-4){1957.317}}
\end{picture}
\end{center}
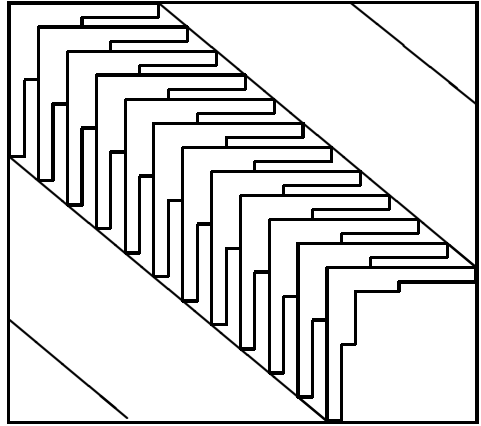
\centerline{Figure 4: An optimal tiling}


\begin{thebibliography}{99}

\bibitem{andrews} G. Andrews, {\bf The Theory of Partitions}, Encyclopedia
 of Mathematics and its Applications, Vol. 2.
Addison-Wesley Publishing Co., Reading, Mass.-London-Amsterdam, 1976.

\bibitem{MS}
F. J. MacWilliams and N. J. A. Sloane,
{\bf The Theory of Error-Correcting Codes},
North Holland, Amsterdam, 1977.

\bibitem{pittel} B. Pittel, On a likely shape of the Random Ferrers diagram,
{\it Adv. Appl. Math}, {\bf 18}, (1997) 432-488. 

\bibitem{szt1} M. Szalay and P. Tur\'an,  On some problems of the statistical 
theory of partitions with application to characters of
the symmetric group. I. Acta Math. Acad. Sci. Hungar. 29 (1977), 
no. 3-4, 361--379. 

\bibitem{szt2} M. Szalay and P. Tur\'an, On some problems of the 
statistical theory of partitions with application to characters of
the symmetric group. II. Acta Math. Acad. Sci. Hungar. 29 (1977),
 no. 3-4, 381--392. 

\bibitem{szt3} M. Szalay and P. Tur\'an, On some problems of the
 statistical theory of partitions with application to characters of
the symmetric group. III. Acta Math. Acad. Sci. Hungar. 32 (1978), 
no. 1-2, 129--155. 

\bibitem{wilf} H. Wilf,  Private communication, 1998.

 \end{thebibliography}
\end{document}